\documentclass[reqno]{amsart}
\usepackage{amssymb}
\usepackage{hyperref}

\begin{document}
\title[ ]
{Spectral equality for $C_0$ semigroups }

\author[ A. TAJMOUATI, M. AMOUCH and M. R. F. ALHOMIDI ZAKARIYA]
{  A. TAJMOUATI, M. AMOUCH and M. R. F. ALHOMIDI ZAKARIYA}  

\address{A. TAJMOUATI and M. R. F. ALHOMIDI ZAKARIYA \newline
 Sidi Mohamed Ben Abdellah
 Univeristy
 Faculty of Sciences Dhar Al Mahraz Fez, Morocco.}
\email{abdelaziztajmouati@yahoo.fr}

\address{Mohamed AMOUCH \newline
Department of Mathematics
University Chouaib Doukkali,
Faculty of Sciences, Eljadida.
24000, Eljadida, Morocco.}
\email{mohamed.amouch@gmail.com}



\newtheorem{theorem}{Theorem}[section]
\newtheorem{definition}{Definition}[section]
\newtheorem{remark}{Remark}
\newtheorem{lemma}{Lemma}[section]
\newtheorem{proposition}{Proposition}[section]
\newtheorem{corollary}{Corollairy}[section]
\newtheorem{example}{Example}
\newcommand{\pr} {{\bf   Proof: \hspace{0.3cm}}}

\maketitle

\begin{abstract}
In this paper,  we give conditions for which the $C_0$ semigroups satisfies
spectral equality for semiregular, essentially semiregular and semi-Fredholm spectrum.
Also, we establish the spectral inclusion for B-Fredholm spectrum of a $C_0$ semigroups.
\end{abstract}

\section{Introduction and preliminaries.}

Let $X$ a Banach space and $\mathcal{B}(X)$ the algebra of all bounded linear operator on $X.$
For $T\in \mathcal{B}(X)$, let $N(T)$,  $R(T)$,  $ R^{\infty}(T)=\bigcap_{n\geq0}R(T^n)$, denote
respectively  the kernel, the range  and the hyper-range of $T.$\\
We denote by $\rho(T)$, $\sigma(T)$, $\sigma_{p}(T)$,  $\sigma_{ap}(T)$, $\sigma_{r}(T)$ and  $\sigma_{su}(T)$
respectively the resolvent set, the spectrum,  the point spectrum, the approximate point spectrum, the residual spectrum and the
surjective spectrum of $T$.\\
For $\lambda\in\rho(T)$, we denote by $R(\lambda,T)=(\lambda -T)^{-1} \in \mathcal{B}(X)$
 the resolvent operator of $T$.\\
 \indent Recall that $T$ is semiregular if $ R(T) $ is closed and $N(T)\subseteq R^{\infty}(T)$ and
 $T$ is essentially semiregular if $ R(T) $ is closed and
 there exists a finite dimensional subspace $G\subseteq X$ such that $N(T)\subseteq R^{\infty}(T)+G.$
 A closed linear operator $T$ is said to be upper semi-Fredholm if $R(T)$ is closed and $dim N(T)<\infty$, and $T$ is the lower semi-Fredholm if $codimR(T)<\infty$, and if the $dim N(T)$  and  $codimR(T)$ are both finite then $T$ is called a Fredholm operator.\\
Let $A\in \mathcal{B}(X)$ and $n$ a nonnegative integer. Define $T_{n}$ to be the restriction of
$T$ to $R(T^{n})$ viewed as a map from $R(T^{n})$ into $R(T^{n})$. If for some integer $n$
the range space $R(T^{n})$ is closed and the induced operator $T_{n}\in \mathcal{B}(R(T^{n}))$
is Fredholm, then $T$ will be said B-Fredholm.
In a similar way, if $T_{n}$ is upper semi-Fredholm (respectively, lower semi-Fredholm )
operator, then $T$ is called upper semi B-Fredholm (respectively, lower semi B-Fredholm),
see \cite{B} for more detail about semi B-Fredholm operators.\\
This classes of operators lead to the definition of the semiregular spectrum, essentially semiregular spectrum,
 semi-Fredholm spectrum, Fredholm spectrum, upper semi-B-Fredholm spectrum,
 lower semi-B-Fredholm spectrum  and the B-Fredholm spectrum defined by:\\
 $$\sigma_{\gamma}(T)=\{\lambda\in\mathbb{C}: T-\lambda I ~~\mbox{is not semiregular}\};$$
 $$\sigma_{\gamma e}(T)=\{\lambda\in\mathbb{C}: T-\lambda I ~~\mbox{is not essentially semiregular }\};$$
 $$\sigma_{\pi}(T)=\{\lambda\in\mathbb{C}: T-\lambda I~~ \mbox{is not semi-Fredholm }\};$$
 $$\sigma_{F}(T)=\{\lambda\in\mathbb{C}: T-\lambda I~~ \mbox{is not Fredholm }\}.$$
 $$\sigma_{uBF}(T)=\{\lambda\in\mathbb{C}: T-\lambda I~~ \mbox{is not  upper semi-B-Fredholm }\};$$
 $$\sigma_{lBF}(T)=\{\lambda\in\mathbb{C}: T-\lambda I~~ \mbox{is not  lower semi-B-Fredholm }\};$$
  $$\sigma_{BF}(T)=\{\lambda\in\mathbb{C}: T-\lambda I~~ \mbox{is not  B-Fredholm }\}.$$
\indent Let $X^{*}$ denote the dual space of $X$ and $A^{*}$ the adjoint operator of  $A$ with domain $D(A)$. Define
the reduced minimum modulus $\gamma(A)$ by setting
$$\gamma(A)=\inf \big \{\frac{\|Au\|}{d(u,N(A))}, ~~u\in D(A)\backslash N(A) \big \}.$$
It is known  that if $D(A)$ is dense in $X,$ then $\gamma(A) = \gamma(A^{*})$ and $\gamma(A) > 0$ if and only if $R(A)$
is closed.\\
Recall that the family $(T(t))_{t\geq0}$ of operators on $X$ is
called a strongly continuous semigroup of operators ( in short $C_0$-semigroups )if:\\
\begin{enumerate}
  \item $T(0)=I;$
  \item $T(s+t)=T(s)T(t)~~~~\mbox{ for all }s,t\geq0;$
  \item $lim_{t\downarrow 0}T(t)x=x~~~~~\mbox{ for every }x\in X.$
\end{enumerate}
Let $\mathcal{T} = (T (t))_{t\geq0}$ be a strongly continuous semigroup on $X.$
The linear operator $A$ defined in \\
$$D(A)=\{x\in X : lim_{t\downarrow 0} \frac{T(t)x-x}{t} \mbox{ exist }\}$$
by\\
$$Ax=lim_{t\downarrow 0}\frac{T(t)x-x}{t}= \frac{d^{+}T(t)x}{dt}\mid_{t=0} \mbox{ for } x\in D(A)$$
is the infinitesimal generator of the semigroup $T(t)$ and $D(A)$ its domain.\\
Let $\mathcal{T} = (T (t))_{t\geq0}$ be a strongly continuous semigroup with generator $A$ on $X.$
We will denote the type (growth bound) of $\mathbb{\tau}$ by
$$\omega_{0}=\lim_{t\rightarrow\infty}\frac{\ln\|T'(t)\|}{t}$$
$$=\inf \big \{\omega \in \mathbb{R}~~: there~~ existe~~ M ~~such ~~that\,\, \|T(t)\|\leq Me^{\omega t },t\geq0 \big \}.$$
\indent Let $A\in \mathcal{B}(X). $ The quasi-nilpotent part and the analytic part play an important
role in local spectral theory, see the monographic of Laursen and  Neumann \cite{LN}.
Recall that the quasi nilpotent part of $A$ denoted
$H_{0}(A)$ is given by
$$H_{0}(A)=\big \{x \in \bigcap_{n\geq 0}D(A^{n})~~:\lim_{n\rightarrow\infty}\|A^{n}x\|^{\frac{1}{n}}=0\big \}.$$
Also, the algebraic core and the analytic core noted respectively by $C(T)$, $K(T)$ are defined by:
$$C(T)=\{x\in X :\exists (x_{n})_{n\geq0}\subset X ,x_{0}=x ~~and ~~Tx_{n}=x_{n-1} \forall n\geq 1\}$$
$$K(T)=\{x\in X :\exists (x_{n})_{n\geq0}\subset X~~and~~ \delta >0 ~~such ~~that~~~ x_{0}=x ~~and$$
$$Tx_{n}=x_{n-1} \forall n\geq 1 ~~and \|x_{n}\|\leq\delta^{n}\|x\| \}.$$
In \cite[Lemma 2.3]{TBK} it is proved that if $A,B\in \mathcal{B}(X) $ are such that $AB = BA$, then :
$$1. K(AB) \subseteq K(A) \bigcap K(B).$$
$$2. C(AB) \subseteq C(T) \bigcap C(B).$$
Recall that some spectral inclusions for various reduced spectra are studied in \cite{En} and \cite{Cha}. The authors proved that
$$e^{t\nu(A)}\subseteq  \nu(T(t))\subseteq e^{t\nu( A)}\bigcup \{0\},$$
where $\nu(.) \in \{\sigma_p(.),\sigma_{ap}(.) ,\sigma_r(.)\}$, point spectrum, approximative spectrum and residual spectrum.\\
\indent In section two of this work, after given the relationship between the analytical core
(respectively, the algebraic core) of a strongly continuous semigroups with its infinitesimal generator,
we will prove the spectral equality for a $C_0$ semigroup
for semiregular, essentially semiregular and semi-Fredholm spectrum, respectively.
In section 3, we will prove the  spectral inclusion  for the B-Fredholm spectrum.\\

\section{Singular spectrum for semigroups generators}

Let $T(t) $ be a $C_{0}$-semigroup on a Banach space $X $ and let $ A$ be its
infinitesimal generator. In this section we will study the relations
between the spectrum of $A $ and the spectrum of each one of the operators
$T(t)$, $t \geq 0$. For this we will begin with the following lemmas proved in \cite{Cha, En, P}.
It will be needed in the sequel.
\begin{lemma}\label{dif}
Let $(A,D(A))$ be the generator of a strongly continuous semigroup
$(T(t))_{t\geq0}$. Then, for every $\lambda \in\mathbb{C}$ and $t > 0$, the following
identities hold:
\begin{enumerate}
 \item  $$(e^{\lambda t}- T(t))x=(\lambda-A)\int_{0}^{t}e^{\lambda(t-s)}T(s)xds ,~~~~~~\lambda \in\mathbb{C},x\in X$$
$$=\int_{0}^{t}e^{\lambda(t-s)}T(s)(\lambda-A)xds  ,~~~~~~\lambda \in\mathbb{C},x\in D(A).$$
  \item  $ R(e^{\lambda t} -T(t))\subseteq R(\lambda-A).$
  \item $N(\lambda-A)=\bigcap_{t\geq 0}N(e^{\lambda t} -T(t))\subseteq N(e^{\lambda t} -T(t)).$
  \item $N(e^{\lambda t} -T(t))=\overline{span}\bigcup_{k\in \mathbb{Z}}N(\lambda+\frac{2\pi i k}{t}-A).$
  \end{enumerate}
  \end{lemma}
\begin{lemma}
Let $(A,D(A))$ be the generator of a strongly continuous semigroup
$(T(t))_{t\geq0}$ and $B(\lambda,t)=\int_{0}^{t}e^{\lambda(t-s)}T(s)xds.$
Then for every $ \lambda \in\mathbb{C}$, $t > 0$ and $n\in \mathbb{N}$,
we have the following:
\begin{enumerate}
   \item $(e^{\lambda t}- T(t))^{n}(x)=(\lambda-A)^{n}B(\lambda,t)^{n}(x) ,~~~~~~\lambda \in\mathbb{C},x\in X$\\
                               $=B(\lambda,t)^{n}(\lambda-A)^{n}(x)  ,~~~~~~\lambda \in\mathbb{C},x\in D(A);$
                                                     \item $ R(e^{\lambda t} -T(t))^{n}\subseteq R(\lambda-A)^{n} ;$
                                                     \item $ N(\lambda-A)^{n}\subseteq N(e^{\lambda t} -T(t))^{n};$
                                                     \item $H_{0}(\lambda-A)\subseteq H_{0}(e^{\lambda t} -T(t)).$
 \end{enumerate}
\end{lemma}

\noindent In this direction; we prove the following lemma.

\begin{lemma}
Let$(A,D(A))$ be the generator of a strongly continuous semigroup
$(T(t))_{t\geq0}$ and $B(\lambda,t)=\int_{0}^{t}e^{\lambda(t-s)}T(s)xds.$
Then, for every$ \lambda \in\mathbb{C}$, $t > 0$ and $n\in \mathbb{N}$,
 we have the following:
\begin{enumerate}
 \item$ K(e^{\lambda t}- T(t))\subseteq K(\lambda-A);$
 \item $ C(e^{\lambda t}- T(t))\subseteq C(\lambda-A).$
                                                   \end{enumerate}
\end{lemma}
\begin{pr}
$(1)$ let $x\in K(e^{\lambda t }-T(t))$, then there exists $(x_{n})_{n\geq 0}\subset X $ and
 $\delta>0$ such that $x_{0}=x$,$(e^{\lambda t }-T(t))x_{n}x_{n-1} $ and $\|x_{n}\|\leq\delta^{n}\|x\|.$\\
Let $(y_{n})_{n\geq0}$ be defined by $$y_{n}=B^{n}(\lambda,t)x_{n},$$
 then $y_{0}=x_{0}=x$,hence $(\lambda-A)y_{n}=(\lambda-A)B^{n}(\lambda,t)x_{n}=(\lambda-A)B(\lambda,t)B^{n-1}(\lambda,t)x_{n}$
as $B^{n}(\lambda,t)x \in D(A^{n})\subseteq D(A),$
then $B^{n}(\lambda,t)\mbox{ and }(\lambda-A)$ commutes for all n and
from $(1),$ we have $(\lambda-A)y_{n}=B^{n-1}(\lambda,t)(\lambda-A)B(\lambda,t)x_{n}=
B^{n-1}(\lambda,t)(e^{\lambda t }-T(t))x_{n}=B^{n-1}(\lambda,t)x_{n-1}=y_{n-1}$
and $$\|y_{n}\|=\|B^{n}(\lambda,t)x_{n} \|\leq \|B^{n}(\lambda,t)\|\delta^{n}\|x\|.$$
Let $\delta^{'}=\|B(\lambda,t)\|\delta $ it follows that $\|y_{n}\|\leq \delta^{n'} \|x\|\Rightarrow x\in K(\lambda-A).$\\
$(2)$ Resulting directly from the $(5)$.
\end{pr}

\noindent In \cite{Et, Te, P} the authors studied and developed a spectral theory for semi groups and their generators.
Precisely, they proved the following two theorems:

\begin{theorem}
For the generator $ A $ of a strongly continuous semigroup
$(T (t))_{t\geq0}$, we have the following spectral inclusion
$$e^{t\nu (A)}\subseteq \nu(T(t)) \setminus \{0\} ;~~\forall~t\geq0,$$
for $\nu(.) \in \{\sigma_{\gamma}(.);\sigma_{\gamma e }(.);\sigma_{\pi}(.);\sigma_{F}(.) \}$\\
\end{theorem}

\noindent and

\begin{theorem}\cite{P}\label{ap}
Let $T( t)$ be a $C_{0}$ semigroup and let $A$ be ~its infinitesimal~generator.Then \\
$$ e^{t\sigma_{ap}(A)}\subset \sigma_{ap}(T(t))\subset e^{t\sigma_{ap}(A)}\bigcup \{0\}.$$\\
More precisely, if $\lambda \in \sigma_{ap}(A) ,$ then $e^{\lambda t} \in \sigma_{ap}(T(t))$ and if $e^{\lambda t} \in \sigma_{ap}(T(t))$ there
exists $k \in \mathbb{Z}$ such that $\lambda_{k} = \lambda + \frac{2\pi ik}{t} \in \sigma_{ap}(A)$.\\
\end{theorem}

\noindent In \cite{P} it is proved for a $C_{0}$-semigroup $T( t)$ with $ A $ its infinitesimal
generator that, if $T(t)$ is differentiable for $t > t_{0}$ and $\lambda \in \sigma(A)$, then
$\lambda e^{\lambda t} \in \sigma(AT(t))$. In the following lemma, we will prove that this result also holds for the point spectrum.

\begin{lemma}
Let $T( t)$ be a $C_{0}$-semigroup and let $ A$ be its infinitesimal
generator. If $T(t)$ is differentiable for $t > t_{0}$ and $\lambda \in \sigma_{p}(A)$, then
$\lambda e^{\lambda t} \in \sigma_{p}(AT(t))$.
\end{lemma}

\begin{pr}
If $\lambda \in \sigma_{p}(A)$, $t > to$, $\exists x \in D(A) , x\neq 0$ such that $(\lambda I- A)x=0.$
Assuming now that $t > t_{0}$ and differentiating the equality 1. of lemma \ref{dif} with respect to $t$,
 we obtain:\\
 $\lambda e^{\lambda t}x - AT(t)x = (\lambda I - A)B^{'}(\lambda ,t)x \mbox{ for every }x \in X.$
$$=B^{'}(\lambda ,t)(\lambda I - A)x \mbox{ for every }x\in D(A).$$
This implies that $\lambda e^{\lambda t}x -(AT(t))x=0,$ that is $\lambda e^{\lambda t} \in \sigma_{p}(AT(t)).$
\end{pr}

\noindent To continue the development of a spectral theory for semi groups and
their generators, we prove that the formula holds for
semiregular spectrum, essentially semiregular spectrum, Fredholm spectrum and semi Fredholm spectrum.
For this we begin with the following result proved in \cite{Vm} which will be used to prove the following theorem:

\begin{proposition}\label{sem}
\begin{enumerate}
  \item If $A\in \mathcal{B}(X)$ is semiregular,then the operator $\widehat{A}:X/R^{\infty}(A)\rightarrow X/R^{\infty}(A)$ induced by $A$ is bounded below.\\
  \item If $A\in \mathcal{B}(X)$ is essentially semiregular, then the operator $\widehat{A}:X/R^{\infty}(A)\rightarrow X/R^{\infty}(A)$ induced by $A$ is
upper semi-Fredholm.
\end{enumerate}
\end{proposition}

\begin{theorem}
For the generator $A $ of a strongly continuous semigroup
$(T (t))_{t\geq0}$, we have
$$ e^{t\nu (A)}\subset \nu(T(t)) \subset e^{t\nu (A)}\bigcup\{0\} ; \mbox{ forall } t\geq 0.$$
More precisely, if $\lambda\in \nu (A),$ then
$e^{\lambda t} \in \nu(T(t))$
and if $e^{\lambda t} \in \nu(T(t))$ and
$(T(t))_{t\geq0}$ is a periodic with period $t,$ then there exists $k \in \mathbb{Z}$ such that
 $\lambda_{k} = \lambda +\frac{2\pi ik}{t} \in \nu(A)$.\\
where $\nu \in \{\sigma_{\gamma};\sigma_{\gamma e };\sigma_{\pi};\sigma_{F}\}$
\end{theorem}
\begin{pr}
if $\lambda\in \nu (A)$then
$e^{\lambda t} \in \nu(T(t))$ see \cite[theorem 2.1]{Et} and \cite[theorem 2]{Sa}.\\
To prove the second inclusion as:\\
For the regular spectrum: Let $t_{0}>0$ be fixed and suppose that for all
$\lambda \in \mathbb{C}\setminus \lambda+2\pi i nt^{-1}$ such that $(\lambda -A)$ is semiregular.
 We show that   $(e^{\lambda t_{0}}-T(t_{0})) $ is semiregular.
 We consider the closed $(T(t))_{t\geq 0}$-invariant subspace $M=R^{\infty}(e^{\lambda t_{0}}-T(t_{0}))$ of $X$
  and quotient semigroup,
$(\widehat{T}(t))_{t\geq 0}$ defined on $X/M$ by:
 $\widehat{T}(t)\widehat{x}=\widehat{T(t)x},$~~~~~~~~~for ~~$\widehat{x}\in X/M,$
  with generator $\widehat{A} $ defined by:
 $$D(A)=\{\widehat{x},x\in D(A)\} ~~ and~~~{A}\widehat{x}=\widehat{Ax}, ~~ for ~~all ~~\widehat{x}\in D(A).$$
 From $(1)$ of proposition$(2.1)$ it flows that the operator $(\lambda -\widehat{A})$ is bounded below for all
  $\lambda \in \mathbb{C}\setminus \lambda+2\pi i kt^{-1}$ and for all $ k \in \mathbb{Z}$.
 Thus, $\lambda \notin \sigma_{ap}(\widehat{A})$.By virtue of Theorem $(2.2)$, we get $e^{\lambda t_{0}}\notin \sigma_{ap}(\widehat{T}(t_{0}))$
 in consequence, the  operator $(e^{\lambda t_{0}}- \widehat{T}(t_{0}))$ is bounded below.
 If $(e^{\lambda t_{0}}-T(t_{0}))x=0$ then $(e^{\lambda t_{0}}- \widehat{T}(t_{0}))(x+M)=0$ and the injective of $(e^{\lambda t_{0}}- \widehat{T}(t_{0}))$
$\Rightarrow x\in M.$ thus $N(e^{\lambda t_{0}}-T(t_{0}))\subset M.$\\
  We found that
 $$ N(e^{\lambda t_{0}}-T(t_{0}))\subset R^{\infty}(e^{\lambda t_{0}}-T(t_{0}))$$
 Now, let us show $R(e^{\lambda t_{0}}-T(t_{0}))$ is closed.To do this, let a sequence$(y_{n})_n$ of elements of $R(e^{\lambda t_{0}}-T(t_{0}))$
 and $y_{n}\rightarrow ~~y$.then there exists a sequence $(v_{n})_{n}\in X$ such that $$R(e^{\lambda t_{0}}-T(t_{0}))v_{n}=y_{n}=(\lambda -A)B(\lambda,t_0)v_{n}=(\lambda -A)u_{n} \rightarrow  y.$$
 As $B(\lambda,t_0) $ is invertible and  $(\lambda -A)$ is closed  then
  $$ (\lambda -A)u_{n}\rightarrow (\lambda -A)u=(\lambda -A)B(\lambda,t_0)B^{-1}(\lambda,t_0)u=(\lambda-A)B(\lambda,t_0)h.$$
 Then $y_{n}=(e^{\lambda t_{0}}-T(t_{0}))v_{n}\rightarrow ~~y=(\lambda-A)B(\lambda,t_0)h=(e^{\lambda t_{0}}-T(t_{0}))h$
 $\Rightarrow y \in R(e^{\lambda t_{0}}-T(t_{0}))$ $\Rightarrow R(e^{\lambda t_{0}}-T(t_{0}))$ is closed,
 consequently the operator $(e^{\lambda t_{0}}-T(t_{0}))$ is semiregular.\\
For the left essential spectrum:
Let $t_{0}>0$ be fixed and suppose that $\lambda \notin \sigma_{\pi}(A)$
for all $\lambda \in \mathbb{C}\setminus \lambda+2\pi i kt^{-1}.$
We show that $ e^{\lambda t_{0}}\notin \sigma_{\pi}(T(t_{0}))$.
By using Lemma$(1.1)(3)$ one has $dim N(\lambda-A)<\infty$ as one has $T(t)$ is periodic we infer that
 $dim N(e^{\lambda t_0}-T(t_0))=dim( \overline{span}\bigcup_{n\in \mathbb{Z}}(\lambda+2\pi i kt_{0}^{-1}-A))<\infty$.\\
For prove that $R(e^{\lambda t_0}-T(t_{0}))$ is closed same prove the first  party.\\
For the essential regular spectrum:  One suppose that $\lambda -A$ is
 essentielly semiregular for all $\lambda \in \mathbb{C}\setminus \lambda+ 2\pi i kt^{-1}$
.We show that $e^{\lambda t_0}-T(t_0)$is essentially semiregular.\\
As $\lambda -A$ is essentielly semiregular from $(2)$ of proposition(2.1), it follows that the operator $\lambda-\widehat{A}$ is upper semi-Fredholm
where $\widehat{A}$ is the generator of the quotient semigroup $(\widehat{T}(t_0))_{t\geq 0}$ defined in the first case.
Thus,$\lambda \notin\sigma_{\pi}(\widehat{A}).$ By virtue of the precedent case we get
 $e^{\lambda t_{0}}\notin \sigma_{\pi}(\widehat{T}(t_{0}))$.Then the operator
 $(e^{\lambda t_{0}}- \widehat{T}(t_{0}))$ is upper semi-Fredholm. Next, let $\pi :X\rightarrow ~~X/M$ be the canonical projection.
By using Lemma$(2.1)(3)$  with $dim(N(e^{\lambda t_{0}}- \widehat{T}(t_{0})))<\infty$, it can be verified that
 $$N(e^{\lambda t_{0}}-T(t_{0}))  \subseteq \pi^{-1}( N(e^{\lambda t_{0}}-\widehat{T}(t_{0})))\subset M+G=R^{\infty}(e^{\lambda t_{0}}- T(t_{0}))+G$$
 for a finite dimensional subspace $G$ of $X$.\\
 In fact,we have $N(e^{\lambda t_{0}}-T(t_{0}))\subset N(e^{\lambda t_{0}}-\widehat{T}(t_{0}))$ and  $\pi( N(e^{\lambda t_{0}}-T(t_{0})))\subseteq  N(e^{\lambda t_{0}}-\widehat{T}(t_{0}))$
 $\Rightarrow N(e^{\lambda t_{0}}-T(t_{0}))\subseteq \pi^{-1} N(e^{\lambda t_{0}}-\widehat{T}(t_{0}))= \pi^{-1} (N(e^{\lambda t_{0}}-T(t_{0}))+M)
 =\pi^{-1} (N(e^{\lambda t_{0}}-T(t_{0})))+M\subset G+M.$\\
The closedness of $R(e^{\lambda t_{0}}- T(t_{0}))$ can be proved in exactly the same way as in the first case.
\end{pr}

\section{B-Fredholm Spectrum for  $C_{0}$-semigroups}

\noindent In \cite{Sa} it is proved the spectral inclusion for semigroups holds for the
Fredholm  spectrum. In this section we will prove that this result also holds for
B-Fredholm spectrum.
Recall from \cite{En} that
If $Y$ is a closed subspace of $X$ such that $T(t)Y\subseteq Y$ for all $t\geq 0$,
that is, if $Y$ is $(T(t))_{t\geq0}$-invariant, then the restrictions
$$T(t)_{|}=T(t)_{|Y}$$ form a strongly continuous semigroup $(T(t)_{|})_{t\geq0}$,
called the subspace semigroup, on the Banach space $Y$.
The part of $A$ in $Y $is the operator $A_|$ defined by $$A_{|}y= Ay$$ with domain
$$D(A_{|})=\{y\in D(A)\bigcap Y: Ay\in Y\}.$$
In other words, $A_{|}$ is the maximal operator induced by $A$ on $Y$ which
coincides with the generator of the semigroup $(T(t)_{|})_{t\geq0 }$ on $Y.$\\
In the sequel, we will prove that the following lemma holds.

\begin{lemma}
Let $A$ the generator of $C_0$ semigroup $T(t)_{t\geq0}$ and $B(\lambda,t)=\int_{0}^{t}e^{\lambda(t-s)}T(s)ds $
is a  linear bounded operator on $X,$ then there exist two operators $C$ and $D $ such that
$(\lambda -A),B(\lambda,t),C,D$ are mutually commuting operators, for all $x\in D(A)$
and $C(\lambda -A)x+DB(\lambda,t)x=Ix$ for $t>0$.
\end{lemma}
\begin{pr}
From lemma \ref{dif} we have:
 $$(e^{\lambda t}- T(t))x=B(\lambda,t)(\lambda-A)x  , \mbox{ for } \lambda \in\mathbb{C} \mbox{ and } x\in D(A).$$
Hence $$(I-e^{-\lambda t} T(t))x=e^{-\lambda t}B(\lambda,t)(\lambda-A)x  ,\mbox{ for } \lambda \in\mathbb{C} \mbox{ and } x\in D(A).$$
Subsequently :
$$Ix= e^{-\lambda t}B(\lambda,t)(\lambda-A)x+e^{-\lambda t} T(t)x  ,\mbox{ for } \lambda \in\mathbb{C} \mbox{ and } x\in D(A).$$
If we take integration for both sides from $0$ to $t$, we find that:\\
$$\int_{0}^{t} Ix ds= \int_{0}^{t} e^{-\lambda s}B(\lambda,s)(\lambda-A)xds+\int_{0}^{t}e^{-\lambda s} T(s)xds  ,~~~~~~\lambda \in\mathbb{C},x\in D(A).$$
Hence, $$Ix=G_{\lambda}(t) (\lambda-A)x+\frac{1}{t}\int_{0}^{t}e^{-\lambda s} T(s)xds  $$
$$=G_{\lambda}(t) (\lambda-A)x+\frac{1}{t}e^{-\lambda t}\int_{0}^{t}e^{\lambda (t-s)} T(s)xds  $$
$$=C (\lambda-A)x+D B(\lambda,t)  ,~~~~~~\lambda \in\mathbb{C},x\in D(A);$$
where $$C=G_{\lambda}(t)=\frac{1}{t}\int_{0}^{t} e^{-\lambda s}B(\lambda,s)ds
\mbox{ and }D=\frac{1}{t}e^{-\lambda t}.$$
From \cite{P}, $ (\lambda -A),B(\lambda,t),C,D$ are mutually commuting operators.
\end{pr}

\begin{lemma}
Let $(\lambda-A)$, $ B(\lambda,t),$ $C$ and $D$ be mutually commuting operators in $D(A)$ such that
$C(\lambda -A)+DB(\lambda,t)=I \mbox{ for } t > 0,$ then
\begin{enumerate}
 \item For every  positive integer $n$ there are $C_n,$ $D_n \in D(A)$ such that\\
  $(\lambda-A)^{n},$ $B^{n}(\lambda,t),$ $C_{n},$ $D_{n}$ are mutually
commuting and $$(\lambda-A)^{n}C_{n} + B^{n}(\lambda,t)D_n = I.$$
                                       \item For every positive integer $n$, we have $$R(e^{\lambda t}-T(t))^{n} = R(\lambda-A)^{n} \bigcap R(B^{n}(\lambda,t))
                                        \mbox{ and }
                                         N((e^{\lambda t}-T(t))^{n}) = N((\lambda-A)^{n}) + N(B^{n}(\lambda,t)).$$
                                        Further more,
                                        $$R^{\infty}(e^{\lambda t}-T(t)) = R^{\infty}(\lambda-A) \bigcap R^{\infty}(B(\lambda,t))
                                       \mbox{and}
                                          N^{\infty}(e^{\lambda t}-T(t)) = N^{\infty}(\lambda-A) + N^{\infty}(B(\lambda,t)).$$
                                       \item $N^{\infty}(\lambda-A) \subset R^{\infty}(B(\lambda,t)) $and $N^{\infty}(B(\lambda,t))\subset  R^{\infty}(\lambda-A)$.
                                     \end{enumerate}
\end{lemma}
\begin{pr}
1) For $\lambda \in\mathbb{C}$ and $x\in D(A)$, we have that:
\begin{eqnarray*} Ix&=&((\lambda-A)C+D B(\lambda,t))^{2n-1}x\\
 &=&\sum_{i=0}^{2n-1} \frac{(2n-1)!}{i!(2n-1-i)!}(\lambda-A)^{i}C^{i}B^{2n-1-i}(\lambda,t))D^{2n-1-i}x\\
&=&(\lambda-A)^{n}C_{n}x + B^{n}(\lambda,t)D_nx.\end{eqnarray*}
This, for some $C_{n},D_{n}$ commuting with $(\lambda-A)^{n},B^{n}(\lambda,t)$.\\
2) From lemma \ref{dif} we have:
$R(e^{\lambda t}-T(t))=R((\lambda-A)B(\lambda,t))\subset R(\lambda-A)\bigcap R(B(\lambda,t))$.
If $y\in R(\lambda-A)\bigcap R(B(\lambda,t))$ with $y= (\lambda-A)x=B(\lambda,t)x^{'}$
  for some $x,x^{'}\in D(A),$ then we set $w = Cx^{'} + Dx$. For this $w$, we have
\begin{eqnarray*} (\lambda-A)w &=& (\lambda-A)Cx^{'} + (\lambda-A)Dx \\ &=& (\lambda-A)Cx^{'} +D (\lambda-A)x \\
 &=& (\lambda-A)Cx^{'} + Dy = (\lambda-A)Cx^{'} + DB(\lambda,t)x^{'} \\ &=& x^{'}.\end{eqnarray*}
 From (1), we have $B(\lambda,t)(\lambda-A)w=B(\lambda,t)x^{'}=y.$
 Hence, we conclude that $R(e^{\lambda t}-T(t))^{n} = R(\lambda-A)^{n} \bigcap R(B^{n}(\lambda,t))$
  and $R^{\infty}(e^{\lambda t}-T(t)) = R^{\infty}(\lambda-A) \bigcap R^{\infty}(B(\lambda,t)).$\\
  Similarly $ N(\lambda-A)+ R(B(\lambda,t))\subset N(e^{\lambda t}-T(t))\bigcap D(A)$.
  If $ x\in N(e^{\lambda t}-T(t))\bigcap D(A)$, then  $(\lambda-A)^{n}C_{n}x + B^{n}(\lambda,t)D_{n}x = Ix$
  where $(\lambda-A)Cx \in N(B(\lambda,t))$ and $B^{n}(\lambda,t)D_{n}x \in N(\lambda -A)$.\\
  Thus, $N(e^{\lambda t}-T(t))\bigcap D(A) = N(\lambda-A) + N(B(\lambda,t)).$ From $(1 ),$ we conclude that
  $N(e^{\lambda t}-T(t))^{n} \bigcap D(A^{n}) = N(\lambda-A)^{n} + N(B^{n}(\lambda,t))$
  and $N^{\infty}(e^{\lambda t}-T(t)) = N^{\infty}(\lambda-A) + N^{\infty}(B(\lambda,t)).$
 If $x\in N(\lambda-A) ,$ then $x=B(\lambda,t)Dx\in R(B(\lambda,t))$, thus $N(\lambda-A)\subset R(B(\lambda,t))$
  and by 1), we have  $N(\lambda-A)^{n}\subset R(B^{n}(\lambda,t))$ for every positive integer $n$. If $m \geq n$ then
  $N(\lambda-A)^{n}\subset N(\lambda-A)^{m} \subset R(B^{m}(\lambda,t))$ so that
  $$N(\lambda-A)^{n} \subset R^{\infty}(B(\lambda,t))
  \mbox{ and }N^{\infty}(\lambda-A) \subset R^{\infty}(B(\lambda,t)).$$
  The inclusion $R^{\infty}(B(\lambda,t))\subset N^{\infty}(\lambda-A)$ follows from the symmetry.
  \end{pr}

\begin{theorem}
For the generator $A$ of a strongly continuous semigroup $(T(t))_{t\geq0}$  one  has the spectral inclusion
$$ e^{t\sigma_{BF}(A)}\subseteq \sigma_{BF}(T(t)).$$
\end{theorem}
\begin{pr}
Suppose that $(e^{\lambda t_0} -T(t_0))$ is B-Fredholm for some $\lambda\in \mathbb{C}\setminus \{0\}$ and $t_{0}>0$, then
there exists a positive integer $n$ such that
 $M= R(e^{\lambda t} -T(t_0))^n$ is closed and the restricted semigroup $(e^{\lambda t} -T(t_0)_{\mid M})$ is Fredholm.\\
We will show that $(\lambda-A)$  is B-Fredholm.
To this end, we show that $R(\lambda-A)^{n}$ is closed.
Let $x\in\overline{R((\lambda-A)^n)}$, then there exist a sequence $u_{k}\in D(A^n)$
such that $(\lambda-A)^{n}u_{k}\rightarrow x,$ hence
$$(e^{\lambda t} -T(t_0))^{n}u_{k} = B_{\lambda}(t_{0})^{n}(\lambda-A)^{n}u_{k} \rightarrow B_{\lambda}(t_{0})^{n}x.$$
Since  $M= R(e^{\lambda t} -T(t_0))^n$ is closed, then
$$(e^{\lambda t} -T(t_0))^{n}u_{k} \rightarrow (e^{\lambda t} -T(t_0))^{n}u= B_{\lambda}(t_{0})^{n}(\lambda-A)^{n}u ,$$
for some $u\in D(A^n)$ this implies that
$$x- (\lambda-A)^{n}u \in N(B_{\lambda}(t_{0}))\subseteq R(\lambda-A)^n$$
so that $x\in R(\lambda-A)^n$, so $R(\lambda-A)^n$ is closed.\\
Now, let us to show that $(\lambda-A_{ \mid R(\lambda-A)^n})$ is Fredholm.
We have that $(e^{\lambda t_0} -T(t_0)_{\mid M})$ is Fredholm with
generator  $(\lambda-A_{\mid M \bigcap D(A)}).$
Since $(e^{\lambda t_0} -T(t_0)_{\mid M})$ is Fredholm, then
$dim N(e^{\lambda t} -T(t_0)_{\mid })< \infty$ implies that
$dim N(e^{\lambda t} -T(t_0))\bigcap R(e^{\lambda t} -T(t_0))^{n}\bigcap D(A))< \infty.$
From \cite[lemma 3.2]{B} and from \cite[Lemma 8]{Aie},  we have:\\
$$N(e^{\lambda t_0}-T(t_0))\bigcap R(e^{\lambda t_0}-T(t_0))^{n} \bigcap D(A) = N(\lambda-A) \bigcap R(\lambda -A)^{n} + N(B(\lambda,t_0))\bigcap R^{n}B(\lambda,t_0),$$ hence
$dim(N(\lambda-A)\bigcap R(\lambda -A)^{n}\bigcap D(A))<dim(N(\lambda-A) \bigcap R(\lambda -A)^{n} + N(B(\lambda,t_0))\bigcap R^{n}(B(\lambda,t_0))\bigcap D(A))$
 $$=dim(N(e^{\lambda t_0}-T(t_0))\bigcap R(e^{\lambda t_0}-T(t_0))^{n} \bigcap D(A))<\infty.$$
Then $dim N(\lambda -A_{\mid R(\lambda-A)^{n}\bigcap D(A)})<\infty.$
Moreover from [10,Lemma 4] \\
$$dim(\frac{R(\lambda-A)^{n}}{R(\lambda-A)^{n+1}})\leq max (\frac{R(\lambda-A)^{n}}{R(\lambda-A)^{n+1}},\frac{R(B^{n}(\lambda,t_0))}{R(B^{n+1}(\lambda,t_0))})
\leq dim(\frac{R(e^{\lambda t_{0}}-T(t_{0}))^{n}}{R(e^{\lambda t_{0}}-T(t_{0}))^{n+1}}).$$
So $\lambda-A$ is B-Fredholm.
\end{pr}

\begin{remark}
By the same argument, we can prove that for the generator $A$ of a strongly continuous semigroup
$(T(t))_{t\geq0}$  one  has the spectral inclusion
$$ e^{t\sigma_{\nu}(A)}\subseteq \sigma_{\nu}(T(t))$$
where $\sigma_{\nu}$ is denote the upper semi-B-Fredholm spectrum, lower semi-B-Fredholm.
\end{remark}

\end{document}